\documentclass[11pt]{article}
\usepackage{amssymb}
\usepackage{epsfig, amsmath, amsfonts}

\def\_{\rule{.3em}{.15ex}}
\newcommand{\ls}[1]
    {\dimen0=\fontdimen6\the\font
     \lineskip=#1\dimen0
     \advance\lineskip.5\fontdimen5\the\font
     \advance\lineskip-\dimen0
     \lineskiplimit=.9\lineskip
     \baselineskip=\lineskip
     \advance\baselineskip\dimen0
     \normallineskip\lineskip
     \normallineskiplimit\lineskiplimit
     \normalbaselineskip\baselineskip
     \ignorespaces
    }
\newtheorem{definition}{Definition}[section]

\newenvironment{proof}{{\em Proof.}}{$\hfill\Box$}

\newtheorem{theorem}[definition]{Theorem}
\newtheorem{corollary}[definition]{Corollary}

\ls{1.25}

\begin{document}

\title{The effect of service time variability on maximum queue lengths in $%
M^X/G/1$ queues}
\author{ Ger Koole $^{\ast}$ \and Misja Nuyens \thanks{%
Department of Mathematics, Vrije Universiteit, De Boelelaan 1081a, 1081 HV
Amsterdam, The Netherlands, \texttt{koole@few.vu.nl} and \texttt{%
mnuyens@few.vu.nl}} \and Rhonda Righter \thanks{%
Department of Industrial Engineering and Operations Research, University of
California, Berkeley, California, USA, \texttt{rrighter@ieor.berkeley.edu} }}
\maketitle

\begin{abstract}
We study the impact of service-time distributions on the distribution of the
maximum queue length during a busy period for the $M^X/G/1$ queue. The
maximum queue length is an important random variable to understand when
designing the buffer size for finite buffer ($M$/$G$/1/$n$) systems. We show
the somewhat surprising result that for three variations of the preemptive
LCFS discipline, the maximum queue length during a busy period is smaller
when service times are more variable (in the convex sense). \vskip 1cm
\end{abstract}

\noindent \textbf{AMS Subject Classification}: Primary 60K25, Secondary
90B22.\newline
\newline
\noindent \textbf{Keywords}: maximum queue length, busy period, service
disciplines, LCFS, variability, stochastic orderings, buffer overflow\newline
\newpage

\section{Introduction}

An important design issue for telecommunication systems and other
applications is determining the buffer size when buffers are finite. We can
better understand the effect of a particular buffer size by understanding
the distribution of the maximum queue length during a busy period in an
infinite buffer system. We give a characterization of the busy-period
maximum queue length, $M$, for the $M^{X}/G/1$ queue for three types of
preemptive LCFS (last-come first-served) disciplines: (i) preempted services
are resumed when service recommences (LCFS-p-resume), (ii) preempted
services must be restarted from scratch when service recommences and a new
service time is chosen from the service-time distribution (LCFS-p-repeat
with resampling), and (iii) preempted services must be restarted from
scratch when service recommences but the total service requirement for a
given customer is the same each time it restarts its service
(LCFS-p-repeat-without-resampling). \ These characterizations of $M$ for
each of the queueing disciplines allow us to show the effect of service-time
distributions on $M,$ as stated in (i)-(iii) below. For a fixed service
discipline, let $M$ and $M^{\prime }$ be the maximum number of customers
during a busy period in two $M^{X}/G/1$ queues with respective generic
service times $S$ and $S^{\prime }$, and with the same arrival rate $\lambda 
$, and the same batch-size distribution. We assume that the distributions of 
$S$ and $S^{\prime }$ are such that the queues are stable. In this paper we
show that the following relations hold; the definitions of the various
stochastic orders can be found in the next section.

\begin{description}
\item[(i)] Under the LCFS-p-resume discipline, if $S^{\prime }\leq _{LT}S$,
then $M^{\prime }\leq _{st}M$.

\item[(ii)] Under the LCFS-p-repeat (with resampling) discipline, if $%
E(e^{-\lambda S^{\prime }})\geq E(e^{-\lambda S})$, then $M^{\prime }\leq
_{st}M$.

\item[(iii)] Under the LCFS-p-repeat-without-resampling discipline, if $%
S^{\prime }\leq _{icv}S$, then $M^{\prime }\leq _{st}M$.
\end{description}

A consequence of our results is the somewhat surprising conclusion that $M$
will be stochastically smaller when service times are more variable (in the
convex sense) under the preemptive LCFS disciplines. Miyazawa (1990) and
Miyazawa and Shanthikumar (1991) show that for the finite-buffer $%
M^{X}/G/1/n $ queue under a non-preemptive discipline, the loss rate, i.e.,
the probability that a random customer is lost, will be larger when service
times are more variable in the convex sense. Our result relates to the loss
rate, but the effect goes in the other direction. That is, we have that for
preemptive LCFS disciplines, $P(M>n)$ is smaller when service times are
larger in the convex sense, where $P(M>n)$ can be interpreted as the
probability of at least one loss during a busy period in the $M^{X}/G/1/n$
queue. See also Chang, Chao, Pinedo, and Shanthikumar (1991).

For other results on the impact of the service time and batch size
distributions on various performance measures of queueing systems, see, for
example, Hordijk (2001), Makowski (1994), and Shanthikumar and Yao (1994),
and the references therein. For other applications of the preempt-repeat
service discipline, see, e.g., Adiri, Frostig, and Rinnooy Kan (1991),
Birge, Frenk, Mittenthal, and Rinnooy Kan (1990), Cai, Sun, and Zhou (2004),
and Cai, Wu, and Zhou (2004).

The paper is organized as follows. We first recall some definitions of
stochastic ordering in the next section. We then study $M$ for each of the
preemptive LCFS disciplines. \ Finally we provide some numerical
illustrations of our results.

\section{Preliminaries}

Recall the following stochastic ordering relations for random variables $X$
and $Y$.

\begin{definition}
$X$ is larger than $Y$ in the stochastic sense, $X\geq _{st}Y$, if $%
E\phi(X)\geq E\phi (Y)$ for all increasing functions $\phi $ for which the
expectations exist.
\end{definition}

Equivalently, $X\geq _{st}Y$ if and only if $P(X>t)\geq P(Y>t)$ for all $t$.

\begin{definition}
$X$ is larger than $Y$ in the convex sense, $X\geq _{cx}Y$, if $E\phi(X)\geq
E\phi (Y)$ for all convex functions $\phi$ for which the expectations exist.
\end{definition}

Note that $X\geq _{cx}Y$ implies $EX=EY$ and $Var(X)\geq Var(Y)$. In this
sense the convex ordering is an ordering of variability in random variables.


\begin{definition}
$X$ is larger than $Y$ in the increasing concave sense, $X\geq _{icv}Y$, if $%
E\phi (X)\geq E\phi (Y)$ for all increasing concave functions $\phi $ for
which the expectations exist.
\end{definition}

\begin{definition}
$X$ is larger than $Y$ in the Laplace-transform sense, $X\geq _{LT}Y$, if $%
E[e^{-\theta X}]\leq E[e^{-\theta Y}]$ for all $\theta >0$ for which the
expectations exist.
\end{definition}

Note that $X\geq _{cx}Y$ implies $X\leq _{icv}Y$, which in turn implies $%
X\leq _{LT}Y$.\newline

Finally, for reasons of brevity we use the following notation. When we say $%
X=[\,Y\,|\,Z=z\,]$, we mean that $P(X=x)=P(Y=x\,|\,Z=z)$ for all $x$.

\section{Preemptive LCFS disciplines}

\subsection{LCFS preempt-resume}

We first consider the $M^{X}/G/1$ queue with the LCFS preempt-resume
(LCFS-p-resume) discipline. That is, the customer that has been in the
system the least amount of time is always served, and newly arriving
customers preempt earlier arrivals already in service. Within a batch
customers are arbitrarily labeled, so that we may think of them as arriving
sequentially, though immediately after each other. Thus, one customer in a
newly arriving batch will be considered the most recent arrival and will
immediately enter service, and the rest of the batch cannot be served until
that customer, as well as all customers arriving in later batches that
preempt that customer, are served.

Customers who resume service after being preempted start their service where
they left off. Hence, a random service with service time $S$ that is
preempted when $t$ units of service have already been received has remaining
service time $[S-t|S>t]$. We also assume service is non-idling. Let $T$ be a
generic interarrival time, where $T$ has an exponential distribution with
rate $\lambda $, and let $X$ be a generic batch size with arbitrary
distribution and mean $\mu $. \ We assume that the queue is stable, $\lambda
\mu ES<1$.

Let customer 0 be the last customer in the first batch in the busy period,
i.e., the first customer to enter service, and let $S_{0}$ be the service
time of customer 0. Let $N=N(S_{0})$ be the number of Poisson batch arrival
times that occur during the service of customer 0, and let $%
N(s)=[N(S_{0})|S_{0}=s]$. Note that the service will be interrupted if $%
N(S_{0})>0$. Let $X_{0}$ be the number of customers in the first batch of
the busy period and define $M(k,n)=[M|X_{0}=k,N=n]$ and $M(k)=[M|X_{0}=k]$,
so that $M(X_{0},N)=M=M(X_{0})$. Let $M_{i}$, $i=1,2,\ldots $, be
i.i.d.~copies of $M$, and define $\max_{i=1,\ldots ,n}M_{i}$ to be $0$ if $%
n=0$. For the LCFS-p-resume discipline we then have the following
characterization of $M(k,n)$.

\begin{theorem}
\label{LCFSeqn} The maximum queue length $M(k,n)$ for the $M^{X}/G/1$ queue
under the LCFS-p-resume discipline satisfies 
\begin{equation}
M(k,n)=_{d}\max \{\max_{i=1,\ldots ,n}M_{i}+k;M(k-1)\},\qquad k\geq 1,\
n\geq 0,
\end{equation}%
where $M(0)=0$, and $M_{i}$, $i=1,2,\ldots $, and $M(k-1)$ are independent.
\end{theorem}

\begin{proof}
We can think of constructing the busy period, conditional on $X_{0}=k$, $%
S_{0}=s$ and $N(s)=n$, as follows. Denote the arrival epochs, on a clock
that only ticks when customer 0 is being served, by $0<t_{1}<\cdots <t_{n}<s$%
. A batch of customers arrives at time $t_{1}$ and starts a new independent
busy period (and stops our clock temporarily), except that there are $k$
more customers in the queue (the original customers) throughout that busy
period. When this first sub-busy period is over, at time $t_{1}+\tau $ say,
then customer 0 returns to service and our clock resumes ticking. Another
batch arrives at time $t_{2}+\tau $, starting a new independent busy period,
and so on, until the $n$ sub-busy periods have completed, as well as the
original service time $s$. Then a new busy period starts with the other $k-1$
customers that arrived in the first batch, and the maximum queue length
during that busy period has the same distribution as $M(k-1)$. Because the
arrival process is memoryless, this construction is stochastically
equivalent to the dynamics of a generic $M^{X}/G/1$ busy period starting
with $k$ customers.
\end{proof}

Let $P(k,b)=P(M(k)\leq b)$, and $P(b)=P(M\leq b)=EP(X_{0},b)$. So $P(0,b)=1$
and $P(0)=0$. 
Using the fact that $E[P(b-k)^{N}]$ is the $z$-transform, or probability
generating function of $N$ evaluated at $z=P(b-k)$, we have from Theorem \ref%
{LCFSeqn} that for $1\leq k\leq b$, 
\begin{eqnarray*}
P(k,b) &=&E[P(M+k\leq b)^{N}]P(k-1,b) \\
&=&E[P(b-k)^{N}]P(k-1,b) \\
&=&Ee^{-\lambda (1-P(b-k))S}P(k-1,b).
\end{eqnarray*}

\begin{corollary}
\label{P(k,b)}For $b\geq k\geq 1$, 
\begin{equation*}
P(k,b)=\prod_{i=1}^{k}Ee^{-\lambda (1-P(b-i))S}.
\end{equation*}%
%
%
%
%
\end{corollary}

If we restrict ourselves to unit batch sizes only, so $X\equiv 1$ and $%
P(b)=P(1,b)$, we have the following corollary.

\begin{corollary}
\label{formule} If $X\equiv 1$, then for $b\geq 1$,%
\begin{equation*}
P(b)=Ee^{-\lambda (1-P(b-1))S}.
\end{equation*}
\end{corollary}

Now we can see how the distribution of $S$ affects $M$.

\begin{theorem}
\label{LCFS-pr} For $M^{X}/G/1$ queues operating under the LCFS-p-resume
discipline, if $S^{\prime }\leq _{LT}S$, then $M^{\prime }\leq _{st}M $. In
particular, if $S^{\prime }\geq _{cx}S$, then $M^{\prime }\leq _{st}M $.
\end{theorem}

\begin{figure}[h]
\begin{center}
\epsfxsize=9cm \epsfbox{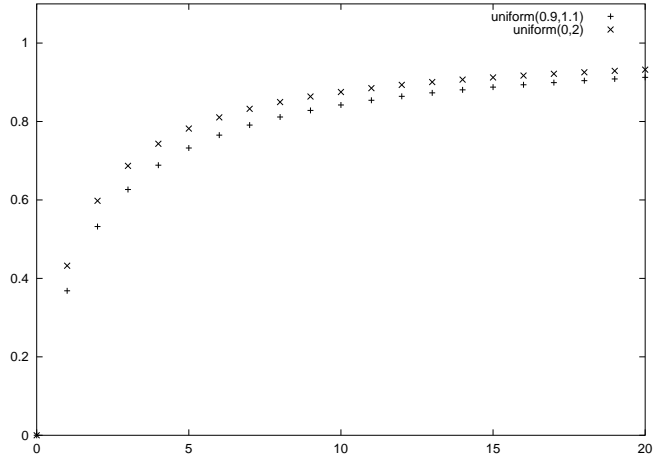}
\end{center}
\caption{$P(M\leq n),n=1,2,\ldots ,20,$ for uniform service-time distributions with $ES=1$%
; $\protect\lambda =0.9.$}\label{figuniform}
\end{figure}

\begin{figure}[h]
\begin{center}
\epsfxsize=9cm \epsfbox{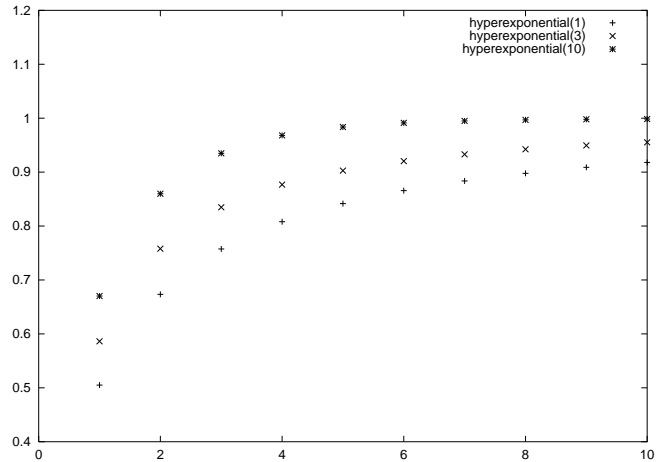}
\end{center}
\caption{$P(M\leq n),n=1,2,\ldots ,10,$ for hyperexponential($k$)
service-time
distributions with $P(S=X_{1})=1-2^{-k}=1-P(S=X_{2})$, where $X_{1}$ and $%
X_{2}$ are exponentially distributed with mean $1/(2(1-2^{-k}))$ and $2^{k}/2
$ respectively, so $ES=1$; $\protect\lambda =0.95$. }
\label{fighyperexp}
\end{figure}

\begin{proof}
To show that $M^{\prime }\leq _{st}M$, we show that $P(k,b)\leq P^{\prime
}(k,b)$ (with the obvious definition for $P^{\prime }$) for all $k$ and $b$
by induction on $b$ and $k$. For each $b$ we have $P(0,b)=1=P^{\prime }(0,b)$%
, and $P(k,1)=0=P^{\prime }(k,1)$ for $k>1$. Since $S^{\prime }\leq _{LT}S$, 
\begin{equation*}
P(1,1)=P(S<T)=Ee^{-\lambda S}\leq Ee^{-\lambda S^{\prime }}=P^{\prime }(1,1).
\end{equation*}%
Suppose $P(i,a)\leq P^{\prime }(i,a)$ for $a<b$ and all $i\geq 0$, so $%
P(a)\leq P^{\prime }(a)$ for all $a<b$, and suppose $P(i,b)\leq P^{\prime
}(i,b)$ for all $0\leq i<k$, and consider $b$ and $k$. From Corollary \ref%
{P(k,b)}, the induction hypothesis, and the assumption $S^{\prime }\leq
_{LT}S$, it then follows that 
\begin{eqnarray*}
P(k,b) &=&\prod_{i=1}^{k}Ee^{-\lambda (1-P(b-i))S}\leq
\prod_{i=1}^{k}Ee^{-\lambda (1-P^{\prime }(b-i))S} \\
&\leq &\prod_{i=1}^{k}Ee^{-\lambda (1-P^{\prime }(b-i))S^{\prime
}}=P^{\prime }(k,b).
\end{eqnarray*}%
This completes the proof.
\end{proof}\newline
\newline
This theorem is illustrated by Figures \ref{figuniform}, \ref%
{figparetos} and \ref{fighyperexp}below, showing
the probabilities $P(M\leq n)$ for (convexly ordered) families of uniform,
Pareto and hyperexponential distributions.\newline

\textbf{Remark } It is well known (Kelly, 1979) that the $M$/$G$/1 queue
under the LCFS-p-resume discipline exhibits service time insensitivity in
the sense that the marginal distribution of the number in the stationary
system, $L$, depends on the service-time distribution only through its mean.
At first this seems at odds with our results, but we must bear in mind that
the maximum number in the system during a busy period depends on the
sample-path evolution of the queue length over a busy period. Hence the
behaviour of $M$ and $L$ may be very different. This idea is further
illustrated by the following heuristic example.\newline

\textbf{Example } Let $M$ be the maximum number in system for an $M$/$G$/1
LCFS-p-resume queue with $S\equiv 1$ (call this system 1) and let $M^{\prime
}$ be the corresponding maximum when the first service time in a busy
period, $S^{\prime }$, is equally likely to be $\varepsilon $ or $%
2-\varepsilon $ so $S\leq _{cx}S^{\prime }$, and the other service times in
the busy period are identically equal to 1 (call this system 2). Then, for $%
\varepsilon $ very small, the first busy period in system 2 is equally
likely to be very short and have a maximum of 1, or it will essentially
consist of two busy periods, each evolving as a busy-period in system 1. The
second of these busy periods starts when the initial customer has received $%
1-\varepsilon/2$ service. That is, roughly, $M^{\prime }$ is equally likely
to be 1 or to have the same distribution as $\max \{M_{1},M_{2}\}$, so $%
M^{\prime }\neq _{st}M$. Note however that $L$ and $L^{\prime }$ have
roughly the same distribution. Indeed, $P(L=0)=P(L^{\prime }=0)$, since the
workload is the same in both systems. Furthermore, a random arrival during a
busy period in system 2 will either see a customer with $S^{\prime
}=\varepsilon $ in service, with very small probability, or will arrive
during one of the two busy periods that each evolve as in system 1. Hence $%
L^{\prime }$ and $L$ have roughly the same distribution. Finally note that
the distribution of the length of a busy period does depend on the
distribution of $S$.\hfill$\Box$

The $M^{X}/G/1/b$ LCFS-p-resume queue also exhibits insensitivity, i.e., the
distribution of the number in system, $L_{b}$, depends on the distribution
of $S$ only through its mean. Hence, the loss rate in the $M^{X}/G/1/b$
queue, $P(L_{b}=b)$, is insensitive to the distribution of $S$. In contrast,
our result shows that the probability of at least one loss during a busy
period, $P(M>b)$, does depend on the distribution of $S$, and is greater
when $S$ is larger in the Laplace-transform sense.

\subsection{LCFS preempt-repeat with resampling}

Now we suppose that when services are preempted they must be restarted from
scratch. The new service time is assumed to be an independent random
variable with the same distribution. We call this the LCFS-p-repeat (with
resampling) discipline. Of course, the behavior of the queue under the
LCFS-p-resume and LCFS-p-repeat disciplines is the same when service times
are exponential.

We use the same notation as in the previous subsection. Now, for stability,
we need $\lambda \mu ES_{e}(S)<1$ and $\lambda \mu ES_{e}(S^{\prime })<1$,
where $S_{e}(S)$ is the effective service time, i.e., the total time a
random customer must spend in service, including restarts due to
interruptions. Thus, 
\begin{equation*}
ES_{e}(S)=E(S\wedge T)+P(S>T)ES_{e}(S),
\end{equation*}%
where $a\wedge b=\min \{a,b\}$, and hence 
\begin{equation}  \label{effective0}
ES_{e}(S)=\frac{E(S\wedge T)}{P(S\leq T)}.
\end{equation}%
For $T$ exponential with rate $\lambda $, it is not hard to show that%
\begin{equation}
ES_{e}=\frac{1-E(e^{-\lambda S})}{\lambda E(e^{-\lambda S})}
\label{effective}
\end{equation}%
and hence for stability we need $E(e^{-\lambda S})>\mu /(\mu +1)$.

For the $M^{X}/G/1$ LCFS-p-repeat queue, we can identify the following
embedded random walk. The number in the system at arrival and departure
epochs during a busy period is equivalent to a random walk on the
nonnegative integers with absorbing state 0. The random walk starts at the
random point $X_{0}$, decreases by 1 if $T>S$ (a departure), and increases
if $T<S$ (an arrival). When it increases, it increases by $X$, where $X$ is
independent of $S$ and $T$. 
Thus, we have the following characterization of $M$, where $I=1$ if $T<S$
and 0 otherwise, and other definitions are as in previous sections.

\begin{theorem}
\label{withoutres} The maximum queue length $M(k)$ for the $M^{X}/G/1$ queue
under the LCFS-p-repeat discipline satisfies 
\begin{equation*}
M(k)=_{d}IM(k+X)+(1-I)\max \{k,M(k-1)\},
\end{equation*}%
where $M(0)=0$, and $I$, $X$, and $M(k-1)$ are mutually independent, and $%
M(k+X)$ is independent of $I$ and $M(k-1)$.
\end{theorem}

Let $I^{\prime }$ be 1 if $T>S^{\prime }$, and 0 otherwise. If $%
P(T>S^{\prime })\geq P(T>S)$, then $I^{\prime }\geq _{st}I$. From Theorem %
\ref{withoutres} and a coupling argument it then follows that $M^{\prime
}\leq _{st}M$. Therefore, we have the following.

\begin{theorem}
\label{LCFS-p-repeat} For $M^{X}/G/1$ queues operating under the
LCFS-p-repeat discipline, if $E(e^{-\lambda S^{\prime}})\geq E(e^{-\lambda
S})$, then $M^{\prime }\leq_{st}M $.
\end{theorem}

Note that for the LCFS-p-repeat discipline, we only need for the Laplace
transform of the service time evaluated at (the arrival rate) $\lambda $ to
be ordered for two service-time distributions, rather than a complete
Laplace-transform ordering. Thus, all possible distributions of service
times can be completely ordered, and hence we have a complete stochastic
ordering of the corresponding maximum queue lengths. Of course, it is also
true that $S^{\prime }\geq _{cv}S$ implies $S^{\prime }\leq _{LT}S$, which
in its turn implies $E(e^{-\lambda S^{\prime }})\geq E(e^{-\lambda S})$.

\subsection{LCFS preemptive repeat without resampling}

For our final model, we suppose again that when services are preempted they
must be restarted from scratch, but now the service time is only drawn from
the service-time distribution once. We call this the
LCFS-p-repeat-without-resampling discipline. Note that the LCFS-p-repeat and
LCFS-p-repeat-without-resampling disciplines are the same for deterministic
service times. For stability, we need again $\lambda \mu ES_{e}(S)<1$ and $%
\lambda \mu ES_{e}(S^{\prime })<1$, where $S_{e}(S)$ is the effective
service time. Given $S=s$, the service time is deterministic and the
effective service time $S_{e}$ is the same as in equation (\ref{effective0}%
), that is 
\begin{equation*}
E[S_e(S)|S=s]=\frac{E(s\wedge T)}{P(s\leq T)}=\frac{1-e^{-\lambda s}}{%
\lambda e^{-\lambda s}}= \frac{1}{\lambda}[e^{\lambda s} -1].
\end{equation*}
Hence, 
\begin{equation*}
ES_e=ES_s(S)=E(E[S_e(S)\,|\, S\,])= \frac{1}{\lambda}[Ee^{\lambda S} -1].
\end{equation*}
So for stability we need $E(e^{\lambda S})<(\mu +1)/\mu $. If, for example, $%
S$ is exponentially distributed with mean $\nu $, then for stability we need 
$\nu >\lambda (\mu +1)$. Note that this value is larger than for the
repeat-with-resampling discipline. Intuitively, a large value of the service
time has a large probability of being interrupted and having to start over,
and each time it restarts it will again have a large service time.

With $X_{0}$, $S_{0}$, and $M$ defined as in the last subsection, and with $%
T_{1}$ defined to be the first interarrival time after the busy period
starts, we now let $M(k,s)=[M|X_{0}=k,S_{0}=s]$ and $M(k)=M(k,S)=[M|X_{0}=k]$%
. Let $I(s)=1$ if $T_{1}<s$ and 0 otherwise. We have the following.

\begin{theorem}
\label{w/o}The maximum queue length $M(k,s)$ for the $M^{X}/G/1$ queue under
the LCFS-p-repeat-without-resampling discipline satisfies 
\begin{equation*}
M(k,s)=_{d}I(s)\max \{M+k;M(k,s)\}+(1-I(s))\max \{k,M(k-1)\},
\end{equation*}%
where $M(0)=0$, and where $I(s)$, $M$, $M(k,s)$, and $M(k-1)$ are
independent.
\end{theorem}

\begin{proof}
Given $X_{0}=k$ and $S_{0}=s$, if an arrival occurs before the first service
completion a new i.i.d.~(sub-)busy period starts, except that there are $k$
additional customers in the queue. When that sub-busy period ends, the
original busy periods starts again,independently of $T_{1}$ and of $M$ for
the ending sub-busy period, with $X_{0}=k$ and $S_{0}=s$. If the first
service completes before an arrival, then we may consider the remainder of
the busy period as a new, independent busy period with $k-1$ initial
customers.
\end{proof}

Let $P(k,b,s)=P(M(k,s)\leq b)$, $P(k,b)=P(M(k)\leq b)=EP(k,b,S_{0})$, and $%
P(b)=P(M<b)=EP(X_{0},b,S_{0})$. We have the following corollary to Theorem %
\ref{w/o}.

\begin{corollary}
For $b\geq k\geq 1$ and for all $s$, 
\begin{equation*}
P(k,b,s)=\frac{P(T>s)P(k-1,b)}{1-P(T<s)P(b-k)}=\frac{e^{-\lambda s}P(k-1,b)}{%
1-(1-e^{-\lambda s})P(b-k)}.
\end{equation*}
\end{corollary}

Using this corollary, we can show the following.

\begin{theorem}
For $M^{X}/G/1$ queues operating under the LCFS-p-repeat-without-resampling
discipline, if $S^{\prime }\leq _{icv}S$, then $M^{\prime }\leq _{st}M$.
Hence, if $S^{\prime }\geq _{cx}S$, then $M^{\prime }\leq _{st}M$.
\end{theorem}

\begin{proof}
From our corollary above, for $b\geq k\geq 1$, and $s\geq 0$, 
\begin{equation*}
P(k,b,s)=\frac{P(k-1,b)}{e^{\lambda s}(1-P(b-k))+P(b-k)}.
\end{equation*}%
It is easy to show that $f(s):=a/(ce^{\lambda s}+d)$ is a decreasing convex
function of $s$ for all $a,c,d\geq 0$, $c+d>0$ (so $-f(s)$ is increasing and
concave). Hence, if $S^{\prime }\leq _{icv}S$ then $-Ef(S^{\prime })\leq
-Ef(S)$ and $Ef(S^{\prime })\geq Ef(S)$. Also note that $P(k,b,s)$ is
increasing in $P(k-1,b)$ and $P(b-k)$ for fixed $s$. The result now follows
using an induction argument similar to the one in the proof of Theorem \ref%
{LCFS-pr}.~
\end{proof}

\section{Numerical Illustrations}

For the LCFS preempt-resume discipline, we calculated $P(M\leq n)$ using
Corollay \ref{formule}. The figures below show the results for several
convexly ordered families of distributions, illustrating Theorem \ref%
{LCFS-pr}.

\begin{figure}[h!]
\begin{center}
\epsfxsize=9cm \epsfbox{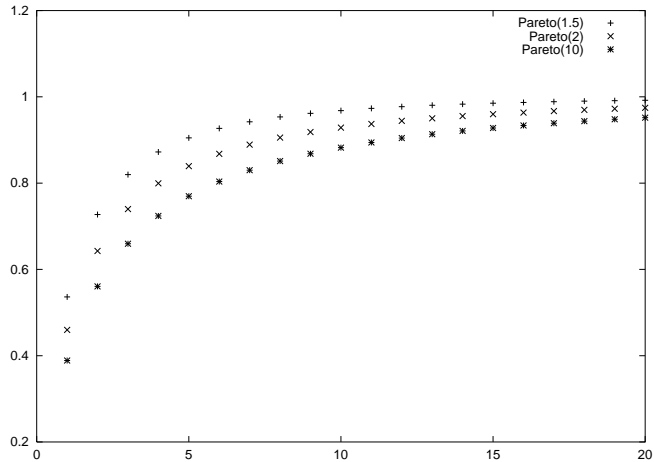}
\end{center}
\caption{$P(M\leq n),n=1,2,\ldots ,20,$ for Pareto($\protect\alpha )$
service-time distributions with distribution function $F_{\protect\alpha }(x)=1-((\protect%
\alpha -1)/(\protect\alpha x))^{\protect\alpha },x\geq (\protect\alpha -1)/%
\protect\alpha $, so $ES=1$; $\protect\lambda =0.95.$}
\label{figparetos}
\end{figure}

\section*{Acknowledgments}
We would like to thank
the referee for a careful review and excellent suggestions.

\end{document}